\documentclass[1p]{elsarticle}
\usepackage{amssymb}
\usepackage{amsfonts}

\newtheorem{theorem}{Theorem}
\newtheorem{conjecture}[theorem]{Conjecture}

\newproof{pf}{Proof}

\begin{document}
\title{Distant sum distinguishing index of graphs}

\author{Jakub Przyby{\l}o\fnref{fn1,fn2}}
\ead{jakubprz@agh.edu.pl, phone: 048-12-617-46-38,  fax: 048-12-617-31-65}

\fntext[fn1]{Financed within the program of the Polish Minister of Science and Higher Education
named ``Iuventus Plus'' in years 2015-2017, project no. IP2014 038873.}
\fntext[fn2]{Partly supported by the Polish Ministry of Science and Higher Education.}

\address{AGH University of Science and Technology, al. A. Mickiewicza 30, 30-059 Krakow, Poland}

\begin{abstract}
Consider a positive integer $r$ and a graph $G=(V,E)$ with maximum degree $\Delta$ and without isolated edges.
The least $k$ so that a proper edge colouring $c:E\to\{1,2,\ldots,k\}$ exists such that $\sum_{e\ni u}c(e)\neq \sum_{e\ni v}c(e)$ for every pair of distinct vertices $u,v$ at distance at most $r$ in $G$ is denoted by $\chi'_{\Sigma,r}(G)$.
For $r=1$ it has been proved that $\chi'_{\Sigma,1}(G)=(1+o(1))\Delta$.
For any $r\geq 2$ in turn an infinite family of graphs is known with $\chi'_{\Sigma,r}(G)=\Omega(\Delta^{r-1})$.
We prove that on the other hand, $\chi'_{\Sigma,r}(G)=O(\Delta^{r-1})$ for $r\geq 2$.
In particular we show that $\chi'_{\Sigma,r}(G)\leq 6\Delta^{r-1}$ if $r\geq 4$.
\end{abstract}


\maketitle

\section{Introduction}

Vertex distinguishing edge colourings have their origins in 
the concept of \emph{irregularity strength}.
This 
graph invariant was designed in~\cite{Chartrand} as a peculiar measure 
of a ``level of irregularity'' of a graph.
A graph or multigraph is called \emph{irregular} itself if all its vertices have pairwise distinct degrees (see~\cite{ChartrandErdosOellermann} for possible alternative definitions).
Note that there are in fact no irregular graphs at all, except the trivial $1$ vertex case.
Thus to capture the degree of 
irregularity of a graph, the authors of~\cite{Chartrand} exploited the fact that there are in turn irregular multigraphs of any order, except order $2$. 
The irregularity strength of a graph $G=(V,E)$, $s(G)$, is then defined as the least $k$ such that 
we are able to construct an irregular multigraph of a given graph by multiplying some of its edges -- each at most $k$ times.
Equivalently, it is the least $k$ so that an edge colouring $c:E\to\{1,2,\ldots,k\}$ exists attributing every vertex $v\in V$ a distinct \emph{weighted degree} defined as:
$$d_c(v):=\sum_{e\ni v}c(e).$$
This shall be also called the \emph{sum at} $v$, see e.g.
\cite{Aigner,Bohman_Kravitz,Lazebnik,Dinitz,Faudree2,Faudree,Frieze,KalKarPf,Lehel,MajerskiPrzybylo2,Nierhoff,Przybylo,irreg_str2}
for exemplary results concerning $s(G)$.
An intriguing local version of the same problem was proposed in~\cite{123KLT}. 
The parameter investigated there differs from $s(G)$ by the reduction of the pairwise distinction requirement only to \emph{adjacent} vertices, and shall be denoted by $s_1(G)$. The well known \emph{1--2--3 Conjecture} presumes that $s_1(G)\leq 3$ for every graph $G$ without isolated edges, see~\cite{123KLT}.
This was investigated e.g. in~\cite{Louigi30,Louigi,123with13}. 
In general it is however thus far only known
that $s_1(G)\leq 5$, see~\cite{KalKarPf_123}. A distance generalization of this problem, introduced in~\cite{Przybylo_distant}
and referring in particular
to the known distant chromatic numbers (see~\cite{DistChrSurvey} for a survey of this topic),
 handles 
a graph invariant $s_r(G)$ (where $r$ is a positive integer), that is the least integer $k$ so that an edge colouring $c:E\to\{1,2,\ldots,k\}$ exists with $d_c(u)\neq d_c(v)$ for every $u,v\in V$ at distance at most $r$ in $G$
, $u\neq v$ --
see also~\cite{Przybylo_distant_edge_probabil}.

The main subject of this paper is the correspondent of $s_r(v)$ 
in the case of \emph{proper} edge colourings.
For any positive integer $r$ and a graph $G=(V,E)$ without isolated edges, by $\chi'_{\Sigma,r}(G)$ we denote the least integer $k$ such that a proper edge colouring $c:E\to\{1,2,\ldots,k\}$ exists with $d_c(u)\neq d_c(v)$ for every $u,v\in V$ with $1\leq d(u,v)\leq r$, where $d(u,v)$ denotes the distance of $u$ and $v$ in $G$. This is called the \emph{$r$-distant sum distinguishing index} of $G$. Such concept 
develops the study on the earlier \emph{neighbour sum distinguishing index} of $G$, $\chi'_\Sigma(G)=\chi'_{\Sigma,1}(G)$, for  which it was conjectured in~\cite{FlandrinMPSW} that $\chi'_{\Sigma}(G) \leq \Delta(G) + 2$ for any
connected graph $G$ of order
at least three different from the cycle $C_5$.
This was asymptotically confirmed in~\cite{Przybylo_asym_optim} and~\cite{Przybylo_asymptotic_note},
where it was showed that
$\chi'_\Sigma(G)\leq (1+o(1))\Delta(G)$, see also~\cite{BonamyPrzybylo,DongWang_mad,FlandrinMPSW,Przybylo_CN_1,Przybylo_CN_2,WangChenWang_planar} for other results concerning $\chi'_\Sigma$.

Exactly the same upper bound as in the case of $\chi'_\Sigma$ above was conjectured to hold 
for the graph invariant $\chi'_a(G)$~\cite{Zhang} (so called \emph{adjacent strong chromatic index} of $G$), i.e. the least $k$ for which a proper edge colouring $c:E\to\{1,2,\ldots,k\}$ exists attributing distinct sets of incident colours to the neighbours in $G$ (see e.g.~\cite{Akbari,BalGLS,Rucinski_regular,Hatami,HocqMont0,HocqMont,Hornak_planar,WangWang,Zhang}
for a number of partial results and upper bounds for this graph invariant, which is one of the most intensively studied subjects within the area), though obviously $\chi'_a(G)\leq \chi'_\Sigma(G)$ for every graph $G$ without isolated edges.
It is however much more challenging to distinguish vertices by sums than by the corresponding sets (even though the conjectured optimal upper bounds are the same in case of the both 
parameters -- $\chi'_\Sigma$ and $\chi'_a$),
what can be easily seen while attempting to apply the probabilistic method. Such approach was e.g. used 
in~\cite{Hatami} to provide an upper bound $\chi'_a(G)\leq \Delta(G)+C$ for all graphs without isolated edges where $C$ is a constant (in particular, if $\Delta(G)$ is large enough, $C=300$ suffices).
In order to 
bring out the fact that 
distinguishing by sums is indeed much more demanding than by sets, one needs to consider distance correspondents of $\chi'_{\Sigma}$ and $\chi'_a$.
It was in particular conjectured in~\cite{Przybylo_distant_Hatami_edge} that for any $r\geq 2$, analogously as in the case of $r=1$,
$\chi'_{a,r}(G) \leq \Delta(G)+C$ 
under minor assumption that $\delta(G)\geq \delta_0$, where $C$ and $\delta_0$ are constants dependent on $r$.
This was confirmed asymptotically and also exactly for some wide graph classes, 
in particular for all regular (and almost regular) graphs with degree large enough, see~\cite{Przybylo_distant_Hatami_edge} for details.
The same certainly does not hold in case of distinguishing by sums, though.
Indeed, from~\cite{Przybylo_distant} follow lower bounds for $\chi'_{\Sigma,r}$ based on 
research concerning so-called Moore bound (see 
e.g. a survey
~\cite{Mirka} concerning this), 
focused on studying the largest possible number of vertices of a graph with maximum degree $\Delta$ and diameter $r$, denoted by $n_{\Delta,r}$.
Namely, it is known that $\chi'_{\Sigma,r}(G)\geq s_r(G)\geq \frac{n_{\Delta,r}}{\Delta}$, hence using e.g. a construction of \emph{undirected de Bruijn graphs} we get for every $r\geq 2$ an infinite family of graph with $\chi'_{\Sigma,r}(G)\geq \Omega(\Delta^{r-1})$, while using an asymptotic result of Bollob\'as and Fernandez de la Vega~\cite{BollobasFernandez} we even obtain for a fixed $\Delta$ an infinite family of graphs with diameter $r$ tending to infinity of order asymptotically equivalent to $\Delta^{r}$ (hence with $\chi'_{\Sigma,r}(G)$ at least asymptotically equivalent to $\Delta^{r-1}$), see~\cite{Przybylo_distant} for details.
Lower bounds of the same form also hold if we narrow our interest down to regular (or almost regular) graphs.
This shows that the difference between the behaviour of $\chi'_{\Sigma,r}$ and $\chi'_{a,r}$ is enormous  
for $r> 2$, what could not be discerned e.g. in case of distinguishing only neighbours (i.e. for $\chi'_{\Sigma}=\chi'_{\Sigma,1}$ and $\chi'_{a}=\chi'_{a,1}$).

In this paper we provide general upper bounds for $\chi'_{\Sigma,r}$ of the same magnitude as the lower ones above.
In particular
we prove that $\chi'_{\Sigma,r}(G)\leq 6\Delta^{r-1}$ for $r\geq 4$ 
and prove the upper bound of order $\Delta^{r-1}$ also in the remaining cases (for $r=2,3$), see Theorem~\ref{mainTheoremDistPropSum} below for details.
These are the first upper bounds of this order for these graph invariants, refining the result from~\cite{Przybylo_distant}, where only slightly better upper bounds (with the same leading ingredient) were proved to hold for the simpler case of non-proper edge colourings, i.e., the graph invariants $s_r(G)$.



\section{Main Result and Proof}

The mentioned above Moore bound, expressing an upper bound for the largest possible number of vertices of a graph with maximum degree $\Delta$ and diameter $r$ is the following (see~\cite{Mirka}):
$$M_{\Delta,r}:=1+\Delta+\Delta(\Delta-1)+\ldots+\Delta(\Delta-1)^{r-1}.$$
Given a graph $G=(V,E)$ with maximum degree $\Delta$ and a vertex $v\in V$, denote by $N^r(v)$ the set of \emph{$r$-neighbours} of $v$, i.e. vertices $u\neq v$ at distance at most $r$ from $v$ in $G$, and note that $d^r(v):=|N^r(v)|\leq M_{\Delta,r}-1\leq \Delta^{r}$ for any $r\geq 1$.
%

\begin{theorem}\label{mainTheoremDistPropSum}
Let $G$ be a graph without isolated edges and with maximum degree
$\Delta\geq 2$, and let $r$ be an integer, $r\geq 2$. 
Then
$$\chi'_{\Sigma,r}(G)\leq 6\left(\frac{M_{\Delta,r}-1}{\Delta}+\Delta-1\right) +\Delta
,$$
hence
$$\chi'_{\Sigma,r}(G)\leq 6\Delta^{r-1}$$
for $r\geq 4$, 
while $\chi'_{\Sigma,3}(G)\leq 6\Delta^2+\Delta$ and $\chi'_{\Sigma,2}(G)\leq 13\Delta-6$. 
\end{theorem}

\begin{pf}
We fix $r\geq 2$ and prove the theorem by induction with respect to the number of vertices of $G$, denoted by $n$.
It is sufficient to show the thesis in the case when $G$ is a connected graph (which is not an isolated edge)
with maximum degree
$\Delta\geq 2$.

For $n=3$ the theorem obviously holds, so assume $n\geq 4$. Denote
$$M=\frac{M_{\Delta,r}-1}{\Delta} ~~~~{\rm and}~~~~ K=M+\Delta-1,$$
and note that then for every $v\in V$ we in particular have
\begin{equation}\label{r-degree-degree-ineq}
d^r(v)\leq d(v)M.
\end{equation}

Suppose first that $G$ contains a vertex $v$ of degree $2$ such that $d(u)\leq 3$ for every $u\in N(v)$.
Let $H=G-v$. By induction hypothesis we may find a desired colouring of every component of $H$ using colours $1,2,\ldots,6K+\Delta$,
where we use colour $1$ on any potential $K_2$-component of $H$ disregarding temporarily a sum conflict between the ends of such $K_2$.
Let $N(v)=\{u_1,u_2\}$. We then greedily choose a colour in $[1,6K+\Delta]$ for the edge $vu_1$ so that the obtained (partial) edge colouring (of $G$) is proper, the (partial) sum at $v$ is distinct from the (partial) sum at $u_2$ and the sum at $u_1$ is distinct from the sums at all its $r$-neighbours in $G$. We are 
able to do this, as the restrictions above
block at most $2+1+3M$ of the available colours.
Then we greedily choose a colour for $vu_2$ from $[1,6K+\Delta]$ so that the obtained edge colouring of $G$ is proper and the sums at $v$ and $u_2$ are distinct from the sums at their respective $r$-neighbours. This is feasible as such restrictions block at most $1+2+2M+3M$ options.
We thus obtain a desired edge colouring of $G$.

Hence we may assume from now on that:
\begin{itemize}
\item[($\ast$)] every vertex of degree $2$ is adjacent with a vertex of degree at least $4$ in $G$.
\end{itemize}
We may also assume that $G$ 
is not a star (as for any $r$, we obviously have $\chi'_{\Sigma,r}(G)=\Delta$ if $G$ is a star). 
Then there are in $G$ two adjacent vertices with degrees at least $2$. We choose a pair of such adjacent vertices that maximizes the sum of their degrees.
By ($\ast$) 
above, the sum of their degrees must equal at least $6$. We set one of these as a root and denote it as $v_n$, and denote the second of these vertices as $v_{n-1}$. Then we continue constructing a spanning tree of $G$ using BFS algorithm, denoting the consecutively chosen vertices by $v_{n-2},v_{n-3},\ldots$, and denote the obtained tree by $F$. This way we also obtain the ordering $v_1,v_2,\ldots,v_{n-1},v_n$ of the vertices of $G$,
such that
$d_G(v_{n-1}),d_G(v_n)\geq 2$, $d_G(v_{n-1})+d_G(v_n)\geq 6$ and every vertex $v_i$ in this sequence except $v_n$ has 
a \emph{forward neighbour} in $G$, i.e. a neighbour $v_j$ of $v_i$ in $G$ with $j>i$.
Analogously we define \emph{backward neighbours} of $v_i$, 
and \emph{forward} and \emph{backward} \emph{$r$-neighbours} of $v_i$ in $G$. 
Moreover, a
\emph{backward} or \emph{forward edge} of $v_i$ shall be any $v_iv_j\in E$ with $j<i$ or $j>i$, resp.,
while by
the \emph{last forward edge} of 
$v_i$ with $i\neq n-1$ we shall mean an edge $v_iv_j\in E$ with the largest $j$.
Note that in fact, due to the use of BFS algorithm, the set of all last forward edges in $G$ equals $E(F)$.
%
%

We first temporarily remove the edges of the spanning tree $F$ of $G$, decreasing the maximum degree of our graph by at least $1$.
Thus by Vizing's Theorem,
we may properly colour the edges of the obtained $G'=G-E(F)$ with integers in $[2K+1,2K+\Delta]$.
Then we assign colour $2K+1$ to all edges of $F$.
The obtained initial edge colouring of $G$ we denote by $c_0$ (note it does not need to be proper due to potential conflicts involving edges in $E(F)$).

We shall be modifying this in order to construct a desired final proper edge colouring
$f:E\to\{1,2,\ldots,6K+\Delta\}$ 
%
in $n-1$ steps, each corresponding to a
consecutive vertex of the sequence $v_1,v_2,\ldots$,
except the last step, within which the weighted degrees of the both $v_{n-1}$ and $v_n$ shall be adjusted. 
From now on the contemporary edge colouring of $G$ shall be denoted by $c$,
hence the contemporary weighted degree of every $v\in V$ shall be denoted by $d_c(v)$,
while by $d(v)$ we shall mean $d_G(v)$.
The moment a given vertex $v_i$ is analyzed (i.e., in step $i$ of the algorithm, or in step $n-1$ in case of $v_n$) we shall associate with it a $2$-element set $D_i=\{s_i,s_i+2K\}$
chosen from the family of pairwise disjoint sets:
$$\mathcal{D}:=\left\{\{b,b+2K\}:b\in\mathbb{N}, 0\leq (b~{\rm mod}~4K)\leq 2K-1\right\}.$$
Ever since $D_i$ is associated with $v_i$ (i.e., before and after steps $i+1,\ldots,n-1$),
we shall require so that $d_c(v_i)\in D_i$.

We shall admit at most two alterations of the colour of every edge $e$ of $G$ except $v_{n-1}v_n$,
whose colour shall be modified only once -- at the end of the construction and via separate rules.
First time, when $e=v_iv_j$ with $i<j$ is a forward edge (of $v_i$, i.e. in step $i$), we shall allow adding to its colour $2K$ (or $0$),
unless $e$ is the last forward edge of $v_i$, when we allow adding to it any integer from the set $\{0,\ldots,2K-1\}$
so that the obtained afterwards $c(e)$ is not congruent to $c(e')$ modulo $2K$ for any adjacent edge $e'$ of $e$ in $G$.
Note that such requirement concerning properness of an edge colouring modulo $2K$ blocks at most $2\Delta-3$ (as every vertex $v_l$ with $l\leq n-1$ has a forward edge) of these available $2K$ options ($2\Delta-2$ instead of $2\Delta-3$ for $v_{n-1}v_n$) -- this leaves at least $2M+1$ options for the colour of any such last forward edge $e$.
Second time, the colour of $e=v_iv_j$ with $i<j$ may be modified when $e$ is a backward edge (of $v_j$, i.e. in step $j$),
when we shall allow only two possible modifications, i.e., adding $2K$ or subtracting $2K$ from its colour (or doing nothing).
Thus the colour of each edge shall always belong to the set
$\{1,\ldots,6K+\Delta\}$.
The main aim of our colour modifications in each (except the last one), say $i$-th, $i\leq n-2$,
step of the algorithm shall be to find a set $D_i=\{s_i,s_i+2K\}\in \mathcal{D}$ disjoint with all the
previously fixed $D_l$ for all $v_l\in N^r(v_i)$ (with $l<i$)
so that we may assure that $d_c(v_i)\in D_i$ via admissible alterations of colours of the edges incident with $v_i$.

Suppose that so far every rule and all our requirements above have been fulfilled, and we are about to analyze $v_i$
(perform $i$-th step of the algorithm), where $i\leq n-2$.
As we have available at least $2M+1$ options non-congruent modulo $2K$ for the colour of the last forward edge of $v_i$, as mentioned above, and a possibility to modify the sum at $v_i$ by exactly $2K$ via the admitted alteration for every of the remaining $d(v_i)-1$ edges incident with $v_i$ in $G$
(indeed, we unconditionally admitted adding 
$2K$ to the colour of a forward edge of $v_i$, except the last one, and may add or subtract $2K$ from the colour of any backward edge $v_lv_i$ of $v_i$ dependent on whether $d_c(v_l)=s_l$ or $d_c(v_l)=s_l+2K$ -- so that $d_c(v_l)$ remained in $D_l$),
we have available at least
\begin{eqnarray}\label{options_for_v_i}
d(v_i)\left(2M+1\right)>2d(v_i)M
\end{eqnarray}
possibilities for $d_c(v_i)$ via admitted alterations of colours of the edges incident with $v_i$.
We have to only make sure that the option that we shall choose out of these does not belong to $D_l$ for any backward $r$-neighbour $v_l$ of $v_i$.
By~(\ref{r-degree-degree-ineq}) this requirement blocks however merely at most $2d(v_i)M$ integers,
hence by~(\ref{options_for_v_i}) we may perform the admissible colour modifications so that afterwards $d_c(v_i)\notin D_l$ for every $v_l\in N^r(v_i)$ with $l<i$.
We then choose $D_i\in \mathcal{D}$ so that $d_c(v_i)\in D_i$. By the definition of $\mathcal{D}$ this guarantees that $D_i$ is disjoint with all $D_l$ such that $v_l\in N^r(v_i)$ and $l<i$.


It is thus sufficient to comment now on the last step of the algorithm within which we simultaneously adjust the sums at $v_{n-1}$ and $v_n$.
We allow to replace the colour of $v_{n-1}v_n$ with any integer in $[1,6K]$ 
which guarantees properness of the obtained edge colouring modulo $2K$.
This requirement itself excludes at most $2\Delta-2$ potential residues 
modulo $2K$ of a colour for $v_{n-1}v_n$,
hence at least $2M$ remain available. 
Let $$R=\{r_1,r_2,\ldots,r_{2M}\}$$ denote a set of exactly $2M$ residues such that
for each $j=1,2,\ldots,2M$, $c(e)\neq r_j~{\rm mod}~{2K}$ for every edge $e$ adjacent with $v_{n-1}v_n$ in $G$.
For the remaining edges (all except $v_{n-1}v_n$) incident with $v_{n-1}$ or $v_n$, which are their backward edges, we similarly as earlier admit
adding or subtracting $2K$ so that afterwards $d_c(v_j)\in D_j$ 
for every $j\leq n-2$. 
(Note that such changes do not influence properness of an edge colouring modulo $2K$.)
It is now enough to prove that the adjustments on the edges incident with $v_{n-1}$ or $v_n$ can be chosen so that
the obtained $d_c(v_{n-1})$ and $d_c(v_n)$ are distinct from the sums at their respective $r$-neighbours.
As vertices $v_l$ at distance at least $2$ from $v_{n-1}$ and $v_n$ shall not change their weighted degrees in this last step,
we may admit $d_c(v_{n-1})\in D_l$ or $d_c(v_n)\in D_l$ this time.
To be strict we shall require that after the last step:
\begin{description}
  \item[(A)] $d_c(v_{n-1})\neq d_c(v_n)$;
  \item[(B)] the sums at $v_{n-1}$ and $v_n$ are distinct from the sums of their respective $r$-neighbours $v_l$ with $l\leq n-2$;
  \item[(C)] neither of the weighted degrees $d_c(v_{n-1})$, $d_c(v_n)$ belongs to any of the sets $D_c(v_l)$ for any $v_l\in N(v_{n-1})\cup N(v_n)$ with $l\leq n-2$.
\end{description}
%
%
%
For each $k\in\{n-1,n\}$, by~(\ref{r-degree-degree-ineq}), the rules (B) and (C) 
may block 
at most
\begin{equation}\label{blocked_integers}
d(v_k)M-1+(d(v_{n-1})-1)+(d(v_n)-1)<(d(v_k)+2)M
\end{equation}
possible weighted degrees 
for $v_k$. Denote the set of these blocked integers (for $d_c(v_k)$) by $I^{(k)}$, and let $J^{(k)}\subset I^{(k)}$ be the subset of these integers in $I^{(k)}$ that would be attainable for $d_c(v_k)$ via admissible modifications of colours of the edges incident with $v_{n-1}$ or $v_n$ (if we disregard rules (A), (B), (C)) and by setting a colour $c(v_{n-1}v_n)\in[1,6K]$ congruent to some residue in $R$.
We then partition this set into 
$2M$ subsets,
$J^{(k)}=J^{(k)}_1\cup J^{(k)}_2\cup\ldots\cup J^{(k)}_{2M}$, where for each $t\in\{1,\ldots,2M\}$, $J^{(k)}_t$ consists
of all these integers from $J^{(k)}$ which could be attained as the weighted degree of $v_k$ only if we used a colour congruent to $r_t$ modulo $2K$ for the edge $v_{n-1}v_n$ (note that there are always $3$ such options 
in the range $[1,6K]$ for $v_{n-1}v_n$).
Set $$j^{(k)}_t: = |J^{(k)}_t|, ~~~~~~a^{(k)}_t: = \frac{j^{(k)}_t}{d(v_k)+2}$$
for $t\in\{1,\ldots,2M\}$ and $k\in\{n-1,n\}$.
Note that there must exist $t'\in\{1,\ldots,2M\}$ such that $a^{(n-1)}_{t'}+a^{(n)}_{t'}<1$. Otherwise,
$$2M\leq \sum_{t=1}^{2M}\left(a^{(n-1)}_{t}+a^{(n)}_{t}\right) = \sum_{t=1}^{2M}a^{(n-1)}_{t}+ \sum_{t=1}^{2M}a^{(n)}_{t},$$
hence at least one of the two sums, say the second one, on the right hand side of the equality above would be at least $M$, but then
$$(d(v_n)+2)M \leq \sum_{t=1}^{2M}a^{(n)}_{t} (d(v_n)+2) = \sum_{t=1}^{2M}|J^{(n)}_t| = |J^{(n)}|\leq |I^{(n)}|,$$
thus we would obtain a contradiction with inequality~(\ref{blocked_integers}).

Now, since $a^{(n-1)}_{t'}+a^{(n)}_{t'}<1$, while $(d(v_{n-1})+2), (d(v_n)+2)\geq 4$ and $d(v_{n-1})+d(v_n)\geq 6$,
then $j^{(n-1)}_{t'}\leq d(v_{n-1})+1$ and $j^{(n)}_{t'}\leq d(v_{n})+1$, and moreover at least one of the following 
must hold:
\begin{description}
  \item[($1^\circ$)] $j^{(n-1)}_{t'} \leq d(v_{n-1})+1$ and $j^{(n)}_{t'}\leq d(v_{n})-2$, or
  \item[($2^\circ$)] $j^{(n-1)}_{t'} = d(v_{n-1})$ and $j^{(n)}_{t'} = d(v_{n})-1$ with $d(v_{n-1})\in \{4,5\}$ and $d(v_n)=2$, or
  \item[($3^\circ$)] $j^{(n-1)}_{t'}\leq d(v_{n-1})-1$ and $j^{(n)}_{t'}\leq d(v_{n})-1$, or
  \item[($4^\circ$)] $j^{(n-1)}_{t'} = d(v_{n-1})-1$ and $j^{(n)}_{t'} = d(v_{n})$ with $d(v_{n-1})=2$ and $d(v_n)\in \{4,5\}$, or
  \item[($5^\circ$)] $j^{(n-1)}_{t'}\leq d(v_{n-1})-2$ and $j^{(n)}_{t'} \leq d(v_{n})+1$.
\end{description}
We shall first try
to fix the final sum at $v_{n-1}$.
For this goal, to colour the edge $v_{n-1}v_n$ we shall use an integer congruent to $r_{t'}$ ($r_{t'}\in [0,2K-1]$) modulo $2K$, i.e. $r_{t'}$, $r_{t'}+2K$ or $r_{t'}+4K$ ($2K$, $4K$ or $6K$ if $r_{t'}=0$).
Such three options, combined with the admitted adjustments for colours of the remaining $d(v_{n-1})-1$ edges incident with $v_{n-1}$
yield $d(v_{n-1})+2$ possible weighted degrees for $v_{n-1}$, which form an arithmetic progression of difference $2K$.

Suppose first that 
($1^\circ$) is true.
Then at least one of these $d(v_{n-1})+2$ possible weighted degrees for $v_{n-1}$, say $d_0$, is not blocked by conditions (B) and (C). Thus we perform admissible modifications of colours of the edges incident with $v_{n-1}$ so that
$d_c(v_{n-1})=d_0$ (fixing $c(v_{n-1}v_n)\equiv r_{t'} ~({\rm mod}~2K)$). As then, via admissible modifications on all edges incident with $v_n$ except $v_{n-1}v_n$ we may generate $d(v_n)$ sums at $v_n$, at most $d(v_{n})-2$ of which might be blocked by (B) and (C) in this case,
we are left with at least two of these, one of which, say $d'_0$ is distinct from $d_0=d_c(v_{n-1})$.
We then perform the admissible alterations of colours of the edges incident with $v_n$ different from $v_{n-1}v_n$ so that $d_c(v_n)=d'_0$ afterwards.
Analogously we proceed in the symmetrical case ($5^\circ$).

Suppose now that 
($2^\circ$) holds.
Then at least two of the possible weighted degrees for $v_{n-1}$, say $d_1$ and $d_2$ with $d_1<d_2$, attainable with $c(v_{n-1})\in \{r_{t'},r_{t'}+2K,r_{t'}+4K\}$ (analogously with $c(v_{n-1})\in \{2K,4K,6K\}$ if $r_{t'}=0$) are not blocked by (B) and (C).
Suppose that for $i=1,2$ we first try to perform any admissible modifications of the edges incident with $v_{n-1}$ so that $d_c(v_{n-1})=d_i$, and then examine $d(v_n)=2$ 
sums attainable at $v_n$ via admissible modifications on 
the only edge incident with it 
and distinct from $v_{n-1}v_n$, denote the set of these by $A_i$, $|A_i|=2$.
Since at most $d(v_n)-1=1$ of these might be blocked due to (B) and (C) in this case, the only possibility which might prevent us from finishing our construction as above is that we have only one available option left for $d_c(v_{n})$ (i.e. exactly 
$1$ is 
blocked by (B) and (C)), and this option is $d_i=d_c(v_{n-1})$.
So suppose it it is the case for $i=1,2$.
Then however, as each $A_i$, $i=1,2$, consists of two elements which differ by exactly $2K$, 
we obtain that 
and 
$d_2=d_1+4K$ and $J_{t'}^{(n)}=\{d_1+2K\}$.
As $d(v_{n-1})\in \{4,5\}$, at least $2$ neighbours of $v_{n-1}$, say $u_1$ and $u_2$, different from $v_n$ are not neighbours of $v_n$ in $G$
(and hence any change of the colour of $v_{n-1}u_1$ or $v_{n-1}u_2$ does not influence the list of sums attainable at $v_n$).
Note now that as $d_2-d_1=4K$, we may assume
that within our examination above the sum $d_1$ for $v_{n-1}$ was attained using the same colours as in the case of the sum $d_2$ on all edges incident with $v_{n-1}$ in $G$ except the edges $v_{n-1}u_1$ and $v_{n-1}u_2$, whose colours had to be bigger by exactly $2K$ in the case of $d_2$, and in particular with the same colour $r^*$ associated to $v_{n-1}v_n$. Suppose then that $r^*-2K\in [1,6K]$ (the reasoning when $r^*+2K\in [1,6K]$ is analogous). Then in the formerly used colouring of the edges incident with $v_{n-1}$ yielding $d_c(v_{n-1})=d_1$ we may introduce two modifications which do not change the sum at $v_{n-1}$, namely we increase the colour of $v_{n-1}u_1$ by $2K$ and decrease the colour of $v_{n-1}v_n$ by $2K$ (to $r^*-2K$).
This way the list of attainable sums at $v_{n}$ via admissible alteration on the only edge incident with $v_n$ and distinct from $v_{n-1}v_n$ shifts from $\{d_1,d_1+2K\}$ to $\{d_1-2K,d_1\}$, hence we may accomodate $d_1-2K\neq d_c(v_{n-1})$ as the sum at $v_{n}$ (since $J_{t'}^{(n)}=\{d_1+2K\}$) and finish the construction of a desired proper edge colouring of $G$.
Again analogously one may proceed in the symmetrical case ($4^\circ$).

Suppose then finally that ($3^\circ$) holds. Then however we have at least $3$ sums attainable at $v_{n-1}$ via admissible colour alterations of the edges incident with $v_{n-1}$ which are not blocked by (B) and (C), denote these by $d_1,d_2,d_3$ with $d_1<d_2<d_3$.
For $i=1,2,3$ one after another let us perform these admissible colour shifts so that $d_c(v_{n-1})=d_i$.
Denote the set of $d(v_n)$ sums attainable at $v_{n}$ via admissible colour alterations of its incident edges other than $v_{n-1}v_n$ by $A_i$ and assume that all of these are blocked
by (B) and (C) except one which
equals exactly $d_i=d_c(v_{n-1})$ (as otherwise, as $j^{(n)}_{t'}\leq d(v_{n})-1$, we may finish constructing our desired colouring by fixing an available not blocked sum different from $d_i$ at $v_n$).
Then however, as each $A_i$ forms an arithmetic progression of difference $2K$ and each contains exactly one of the three sums $d_1<d_2<d_3$, not blocked for $v_n$ (i.e. $d_1\in A_1, d_2\in A_2, d_3\in A_3$), 
we obtain that $A_1\cap A_3=\emptyset$, hence $|A_1\cup A_2\cup A_3|\geq 2d(v_n)+1\geq d(v_n)+3$, and thus at least one $A_i$ must contain at least $2$ sums not blocked for $v_n$ by (B) and (C), as $j^{(n)}_{t'}\leq d(v_{n})-1$ in this case, what yields a contradiction with our assumption above concerning $A_1,A_2,A_3$.

At the end of our construction we set $f(e)=c(e)$ for $e\in E$ to obtain a desired final proper edge colouring of $G$.
$\blacksquare$
\end{pf}

\section{Conclusion}

We conclude the paper by posing the following two conjectures.
\begin{conjecture}
For every integer $r\geq 3$ and each graph $G$ without isolated edges of maximum degree $\Delta$,
$\chi'_{\Sigma,r}(G)\leq (1+o(1))\Delta^{r-1}$. 
\end{conjecture}
\begin{conjecture}
For every graph $G$ without isolated edges of maximum degree $\Delta$, $\chi'_{\Sigma,2}(G)\leq (2+o(1))\Delta$.
\end{conjecture}

\end{document}